\documentclass[12pt]{article}

\usepackage{lineno,hyperref}
\modulolinenumbers[5]

\usepackage{caption}
\usepackage{amsmath}
\usepackage{amssymb}
\usepackage[numbers]{natbib}
\usepackage{bbm, dsfont}
\usepackage{algorithm}
\usepackage{graphicx}
\usepackage{epstopdf}
\usepackage{tabularx}
\usepackage{rotating}
\usepackage{multirow}
\usepackage{color, colortbl}
\usepackage{algpseudocode}
\usepackage{float}

\usepackage{caption}
\DeclareCaptionType{copyrightbox}
\usepackage{subcaption}
\usepackage{epstopdf}

\newcommand{\ub}{{\bf u}}

\newcommand{\T}{{\bf T}}
\newcommand{\QEST}{{\bf \widehat{\bf Q}}}
\newcommand{\s}{{\bf s}}

\newcommand{\ra}{\zeta}
\newcommand{\zeros}{{\bf 0}}
\newcommand{\BE}{\widehat{\bf B}}

\newcommand{\Nens}{N} 
\newcommand{\Nobs}{m} 
\newcommand{\Nstate}{n} 
\newcommand{\X}{{\bf X}} 
\newcommand{\x}{{\bf x}} 
\newcommand{\lp}{\left (} 
\newcommand{\rp}{\right )} 
\newcommand{\lb}{\left [} 
\newcommand{\rb}{\right ]} 

\newcommand{\B}{{\bf B}} 
\newcommand{\R}{{\bf R}} 
\newcommand{\N}{M} 
\newcommand{\y}{{\bf y}} 
\renewcommand{\H}{{\bf H}} 
\newcommand{\xm}{{\overline{\bf x}}} 
\newcommand{\w}{{\bf r}} 
\newcommand{\I}{{\bf I}} 
\newcommand{\M}{\mathcal{M}} 
\newcommand{\Nor}{\mathcal{N}} 
\newcommand{\Y}{{\bf Y}} 
\newcommand{\lle}{\left \{ } 
\newcommand{\rle}{\right \}} 

\newcommand{\DX}{{{\boldsymbol \delta} {\bf X}}}
\newcommand{\Q}{{\bf Q}} 

\renewcommand{\P}{{\bf P}} 
\renewcommand{\Re}{\mathbbm{R}}
\renewcommand{\S}{{\bf S}} 
\newcommand{\errobs}{{\boldsymbol \epsilon}} 

\newcommand{\D}{{\bf D}}

\newcommand{\Z}{{\bf Z}} 

\newcolumntype{N}{>{\centering\arraybackslash} m{0.30\textwidth} }
\newcolumntype{V}{>{\centering\arraybackslash} m{0.02\textwidth} }

\newcommand{\zero}{{\bf 0}}

\newcommand{\BEST}{\widehat{\bf B}}

\newcommand{\U}{{\bf U}}
\newcommand{\V}{{\bf V}}

\newcommand{\ones}{{\bf 1}}

\title{A Parallel Implementation of the Ensemble Kalman Filter Based on Modified
Cholesky Decomposition}
\author{
        Elias D. Nino \\
        Computational Science Laboratory\\
        Department of Computer Science\\
        Virginia Tech \\
        Blacksburg, VA 24060, \underline{USA}, E-mail: {enino@vt.edu}
            \and
        Adrian Sandu \\
        Computational Science Laboratory\\
        Department of Computer Science\\
        Virginia Tech \\
        Blacksburg, VA 24060, \underline{USA}, E-mail: {asandu7@vt.edu}
         \and
        Xinwei Deng \\
        Department of Statistics \\
        Virginia Tech \\
        Blacksburg, VA 24060, \underline{USA}, E-mail: {xdeng@vt.edu}
}
\date{\today}

\begin{document}

\thispagestyle{empty}
\setcounter{page}{0}

\begin{Huge}
\begin{center}
{\bf Computer Science Technical Report CSTR-3/2016 }\\
\today
\end{center}
\end{Huge}
\vfil
\begin{huge}
\begin{center}
Elias D. Ni\~no, Adrian Sandu and Xinwei Deng
\end{center}
\end{huge}

\vfil
\begin{huge}
\begin{it}
\begin{center}
A Parallel Implementation of the Ensemble Kalman Filter Based on Modified
Cholesky Decomposition
\end{center}
\end{it}
\end{huge}
\vfil

\begin{large}
\begin{center}
Computational Science Laboratory \\
Computer Science Department \\
Virginia Polytechnic Institute and State University \\
Blacksburg, VA 24060 \\
Phone: (540)-231-2193 \\
Fax: (540)-231-6075 \\ 
Email: \url{sandu@cs.vt.edu} \\
Web: \url{http://csl.cs.vt.edu}
\end{center}
\end{large}

\vspace*{1cm}

\begin{tabular}{ccc}
\includegraphics[width=2.5in]{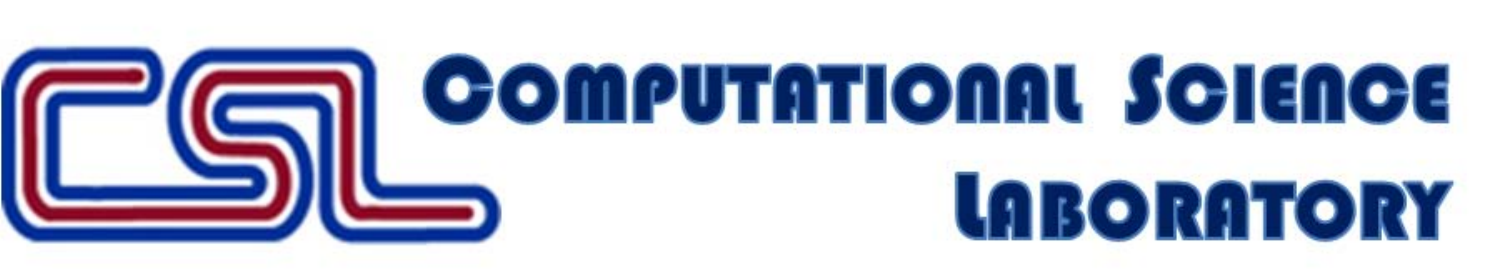}
&\hspace{2.5in}&
\includegraphics[width=2.5in]{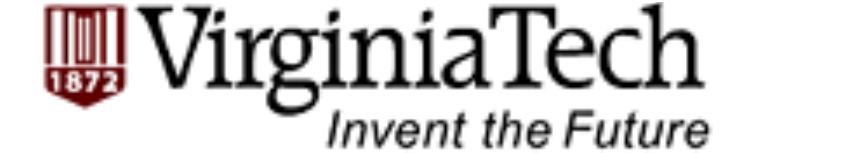} \\
{\bf\em Innovative Computational Solutions} &&\\
\end{tabular}

\newpage

\maketitle

\begin{abstract}
This paper discusses an efficient parallel implementation of the ensemble Kalman filter based on the modified Cholesky decomposition. The proposed implementation starts with decomposing the domain into sub-domains. In each sub-domain a sparse estimation of the inverse background error covariance matrix is computed via a modified Cholesky decomposition; the estimates are computed concurrently on separate processors. The sparsity of this estimator is dictated by the conditional independence of model components for some radius of influence. Then, the assimilation step is carried out in parallel without the need of inter-processor  communication. Once the local analysis states are computed, the analysis sub-domains are mapped back onto the global domain to obtain the analysis ensemble. Computational experiments are performed using the Atmospheric General Circulation Model (SPEEDY) with the T-63 resolution on the Blueridge  cluster at Virginia Tech. The number of processors used in the experiments ranges from 96 to 2,048. The proposed implementation outperforms in terms of accuracy the well-known local ensemble transform Kalman filter (LETKF) for all the model variables. The computational time of the proposed implementation is similar to that of the parallel LETKF method (where no covariance estimation is performed). Finally, for the largest number of processors, the proposed parallel implementation is 400 times faster than the serial version of the proposed method.
\end{abstract}


\section{Introduction}

In operational data assimilation, sequential and variational methods are required to posses the ability of being performed in parallel \cite{nerger2013software,rao2016time,liu2012advancing}. This obeys to current atmospheric and oceanic model resolutions in which the total number of components arises to the order of millions and the daily information to be assimilated in the order of terabytes \cite{DAIssues,emerick2013ensemble}. Thus, serial data assimilation methods are impractical under realistic operational scenarios. In sequential data assimilation, one of the best parallel ensemble Kalman filter (EnKF) implementations is the local ensemble transform Kalman filter (LETKF) \cite{TELA:TELA076}. This method is based on domain localization given a radius of influence $\ra$. Usually, the assimilation process is performed for each model component in parallel making use of a deterministic formulation of the EnKF in the ensemble space. In this formulation, the unknown background error covariance matrix is estimated by the rank-deficient ensemble covariance matrix which, in ensemble space, is well-defined. The LETKF relies in the assumption that local domain analyses avoid the impact of spurious correlations, for instance, by considering only small values for $\ra$. However, in operational data assimilation, $\ra$ can be large owing to circumstances such as sparse observational networks and/or long distance data error correlations (i.e., pressure fields) In such cases, the accuracy of the LETKF can be negatively impacted owing to spurious correlations. 

We think there is an opportunity to provide a more robust parallel ensemble Kalman filter implementation via a better estimation of background error correlations. When two model components (i.e., grid points) are assumed to be conditionally independent, their corresponding entry in the estimated inverse background error covariance matrix is zero. Conditionally dependence/independence of model components can be forced making use of local domain analyses. For instance, when the distance of two model components in physical space is larger than $\ra$, their corresponding entry in the inverse background error covariance matrix is zero. This can be exploited in order to obtain sparse estimators of such matrix which implies huge savings in terms of memory and computations. Even more, high performance computing can be used in order to speedup the assimilation process: the global domain can be decomposed according to an available number of processors, for all processors, local inverse background error covariance matrices are estimated and then,  the stochastic EnKF formulation \cite{EnKFEvensen} can be used in order to compute local domain analyses. The local analyses and then mapped back onto the global domain from which the global analysis state is obtained.

This paper is organized as follows. In section \ref{sec:preliminaries} basic concepts regarding sequential data assimilation and covariance matrix estimation are presented, in section \ref{sec:parallel-implementation}  a parallel implementation of the ensemble Kalman filter based on the modified Cholesky decomposition is proposed; experimental results are discussed in section \ref{sec:experiments} and future research directions are presented in section \ref{sec:future-work}. Conclusions are drawn in section \ref{sec:conclusions}.

\section{Preliminaries}
\label{sec:preliminaries}
\subsection{Modified Cholesky decompositon}
\label{subsec:modified-cholesky}

Let $\S = \lle \s_1,\, \s_2,\,\ldots,\,\s_{\Nens} \rle 
\in \Re^{\Nstate \times \Nens}$, the matrix whose columns are $\Nstate$-th dimensional random Gaussian vectors with probability distribution $\Nor \lp \zero_{\Nstate},\, \Q \rp$, where the number of columns $\Nens$ denotes the number of samples. Denote by $\x^{[j]} \in \Re^{\Nens \times 1}$, the vector holding the $j$-th component across all the columns of $\S$, for $2 \le j \le \Nstate$. The modified Cholesky decomposition \cite{modifiedCholesky} arises from regressing each variable $\x^{[j]}$ on its predecessors $\x^{[j-1]}$, $\x^{[j-2]}$, $\ldots$, $\x^{[1]}$, that is , fitting regressions:
\begin{eqnarray}
\label{eq:modified-Cholesky-decomposition}
\displaystyle
\x^{[j]} = \sum_{q=1}^{j-1} \beta_{jq} \cdot \x^{[q]} + \varepsilon^{[j]} \in \Re^{\Nens \times 1} ,\,
\end{eqnarray}
where $\varepsilon^{[j]}$ denotes the error in the regression of the $j$-th component. Let $\D_{jj} = \lle {\bf var} \lp \varepsilon^{[j]} \rp \rle \in \Re^{\Nstate \times \Nstate}$ be the diagonal matrix of error variances and let $\T_{jq} = \{-\beta_{jq}\} \in \Re^{\Nstate \times \Nstate}$ denote the unitary lower-triangular matrix containing the negative value of  regression coefficients, for $2 \le q < j \le \Nstate$.
 An approximation of the inverse covariance matrix $\Q^{-1} \in \Re^{\Nstate \times \Nstate}$ reads:
\begin{eqnarray}
\label{eq:approximated-Q-inverse}
\Q^{-1} \approx \QEST^{-1} = \T^T \cdot \D^{-1} \cdot \T  \,,
\end{eqnarray}
and making use of basic linear algebra, an approximation of $\Q \in \Re^{\Nstate \times \Nstate}$ is:
\begin{eqnarray}
\label{eq:approximated-Q}
\Q \approx \QEST = \T^{-1} \cdot \D \cdot \T^{-T}  \,.
\end{eqnarray}
\subsection{Local ensemble transform Kalman filter}
\label{subsec:local-ensemble-transform}

\begin{figure}[htp]
\centering
\begin{subfigure}[b]{0.45\textwidth}
\centering
\includegraphics[width=0.7\textwidth]{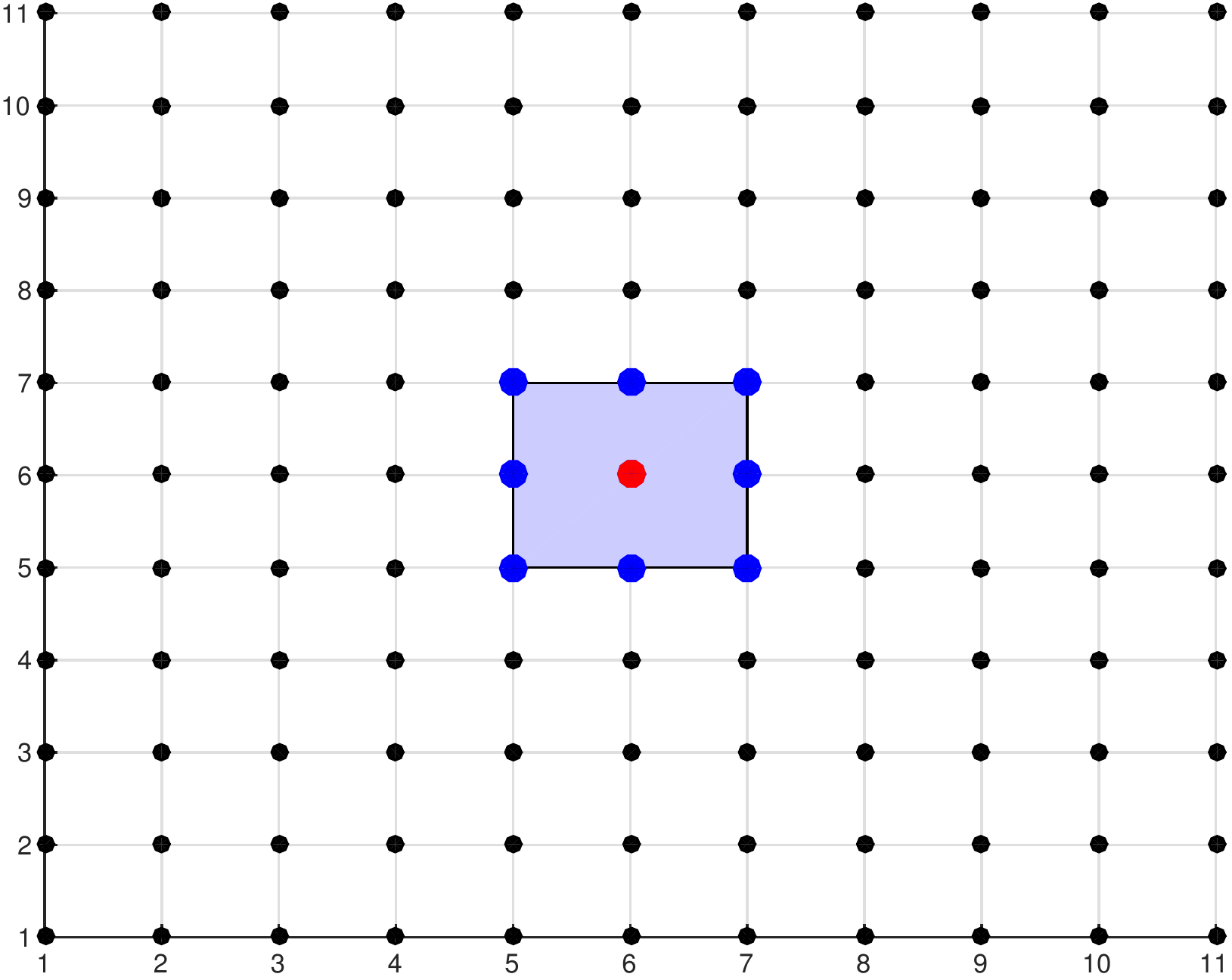}
\caption{$\ra=1$}
\end{subfigure}\hspace{1em}
\begin{subfigure}[b]{0.45\textwidth}
\centering
\includegraphics[width=0.7\textwidth]{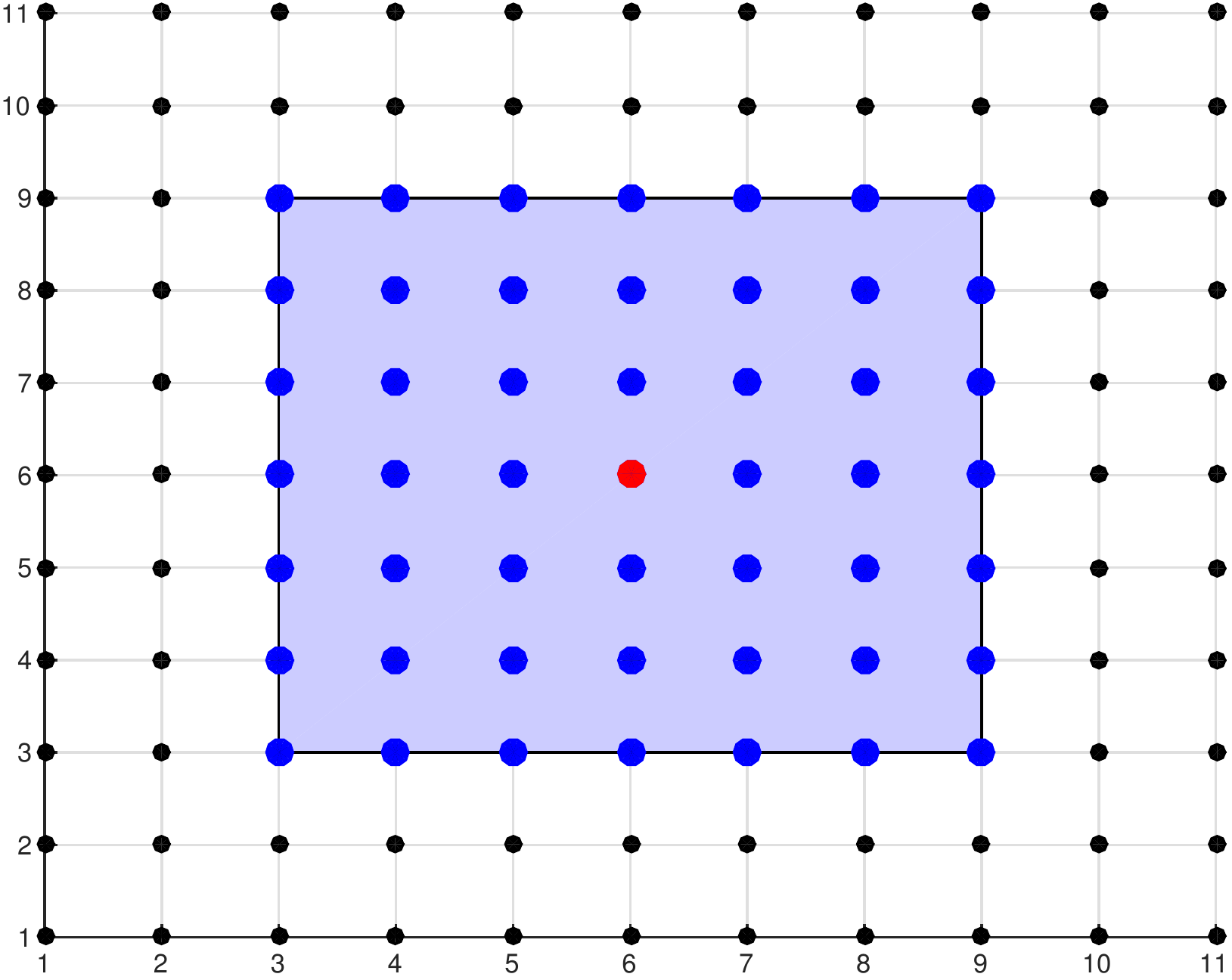}
\caption{$\ra=3$}
\end{subfigure}
\caption{Local boxes for different radius of influence $\ra$.}
\label{fig:local-boxes}
\end{figure}

Localization is commonly used in the context of sequential data assimilation in order to mitigate the impact of spurious correlations in the assimilation process. In general, two forms of localization methods are used: covariance matrix localization and domain localization, both have proven to be equivalent \cite{localization_methods}. In practice, covariance matrix localization can be very difficult owing to the explicit representation in memory of the ensemble covariance matrix. On the other hand, domain localization methods avoid spurious correlations by considering only observations within a given radius of influence $\ra$: in the two-dimensional case, each model component is surrounded by a local box of dimension $(2 \cdot \ra+1,\,2\cdot \ra+1)$ and the information within the scope of $\ra$ (observed components and background error correlations) is used in the assimilation process and conversely, the information out the local box is discarded. In figure \ref{fig:local-boxes}, local boxes for different radii of influence $\ra$ are shown. The red grid point is the one to be assimilated, blue points are used in the assimilation process while black points are discarded. Based on this idea, the local ensemble transform Kalman filter is proposed (LETKF) \cite{TELA:TELA076} 
The global formulation of the LETKF is defined as follows: for a given background ensemble
\begin{eqnarray}
\label{eq:background-ensemble}
\displaystyle
\X^b = \lb \x^{b[1]},\, \x^{b[2]},\, \ldots,\, \x^{b[\Nens]} \rb \in \Re^{\Nstate \times \Nens} ,\,
\end{eqnarray}
and ensemble perturbation matrix
\begin{eqnarray}
\label{eq:ensemble-perturbations}
\displaystyle
\U^b = \X^b - \xm^b \otimes \ones_{\Nens}^T \in \Re^{\Nstate \times \Nens} ,\,
\end{eqnarray}
where $\Nstate$ is the number of model components, $\Nens$ is the ensemble size, $\x^{b[i]} \in \Re^{\Nstate \times 1}$ is the $i$-th ensemble member, for $1 \le i \le \Nens$, $\xm^b$ is the ensemble mean, $\ones_{\Nens}$ is the $\Nens$-th dimensional vector whose components are all ones and $\otimes$ denotes the outer product of two vectors, an estimated of the analysis error covariance matrix in the ensemble space reads:
\begin{subequations}
\label{eq:LETKF-analysis}
\begin{eqnarray}
\widehat{\P^a} = \lb \lp \Nens-1\rp \cdot \I_{\Nens \times \Nens} + \Z^T \cdot \R^{-1} \cdot \Z \rb^{-1}
\end{eqnarray}
where $\Z = \H \cdot \U^b \in \Re^{\Nobs \times \Nens}$, $\H \in \Re^{\Nobs \times \Nstate}$ is the linear observational operator, $\Nobs$ is the number of observed components and, $\R \in \Re^{\Nobs \times \Nobs}$ is the estimated data error covariance matrix. The optimal weights in such space reads:
\begin{eqnarray}
\displaystyle
\w^a = \widehat{\P^a} \cdot \Z^T \cdot \R^{-1} \cdot \lb \y - \H \cdot \xm^b \rb,\,
\end{eqnarray}
therefore, the optimal perturbations can be computed as follows:
\begin{eqnarray}
\displaystyle
{\bf W}^a = \w^a \otimes \ones_{\Nens}^T + \lb (\Nens-1) \cdot \widehat{\P^a} \rb^{1/2} \in \Re^{\Nens \times \Nens}
\end{eqnarray}
from which, in model space, the analysis reads:
\begin{eqnarray}
\label{eq:analysis-LETKF}
\displaystyle
\X^a = \xm^b \otimes \ones_{\Nens}^{T} + \U \cdot {\bf W}^a \in \Re^{\Nstate \times \Nens} \,.
\end{eqnarray}
\end{subequations}
The set of equations \eqref{eq:LETKF-analysis} are applied to each model component in order to compute the global analysis state.

\subsection{Ensemble Kalman Filter Based On Modified Cholesky}
\label{sec:proposed-method}

\begin{figure}[H]
\centering
\begin{subfigure}{0.4\textwidth}
\centering
\includegraphics[width=0.6\textwidth,height=0.6\textwidth]{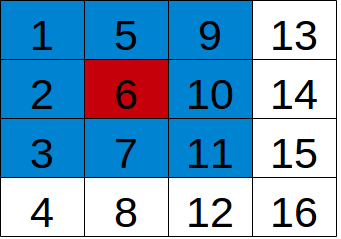}
\caption{In blue, local box for the model component 6 when $\ra=1$.}
\end{subfigure} \hspace{1em}
\begin{subfigure}{0.4\textwidth}
\centering
\includegraphics[width=0.6\textwidth,height=0.6\textwidth]{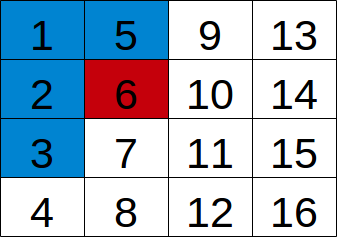}
\caption{In blue, predecessors of the model component 6 for $\ra=1$.}
\end{subfigure}
\caption{Local model components (local box) and local predecessors for the model component 6 when $\ra=1$. Column-major ordering is utilized to label the model components.}
\label{fig:predecessors}
\end{figure}

In \cite{nino2016enkf-mc}, the modified Cholesky decomposition is used in order to obtain sparse estimators of the inverse background error covariance matrix. the columns of matrix \eqref{eq:ensemble-perturbations} are assumed normally distributed with moments:
\begin{eqnarray}
\ub^{b[i]} \sim \Nor \lp \zero_{\Nstate},\, \B \rp, \,\text{ for $1 \le i \le \Nens$},\,
\end{eqnarray}
where $\B \in \Re^{\Nstate \times \Nstate}$ is the true unknown background error covariance matrix. Denote by $\x^{[j]} \in \Re^{\Nens \times 1}$ the vector holding the $j$-th model component across all the columns of matrix \eqref{eq:ensemble-perturbations}, for $1 \le j \le \Nstate$, following the analysis of section \ref{subsec:modified-cholesky}, i.e., $\S = \U$, an estimate of the inverse background error covariance matrix reads:
\begin{eqnarray}
\label{eq:inverse-B}
\displaystyle
\B^{-1} \approx \BEST^{-1} =  \T^T \cdot \D^{-1} \cdot \T \in \Re^{\Nstate \times \Nstate},\,
\end{eqnarray}
and similar to \eqref{eq:approximated-Q}, 
\begin{eqnarray}
\label{eq:B}
\displaystyle
\B \approx \BEST = \T^{-1} \cdot \D \cdot \T^{-T} \in \Re^{\Nstate \times \Nstate} \,.
\end{eqnarray}
Based on \eqref{eq:modified-Cholesky-decomposition}, the resulting estimator $\BEST^{-1}$ can be dense. This implies no conditional independence of model components in space which, in practice, can be quite unrealistic for model variables such as wind components, specific humidity and temperature. Thus, a more realistic approximation of $\B^{-1}$ implies a sparse estimator $\BEST^{-1}$. Readily, the structure of $\BEST^{-1}$ depends on the structure of $\T$ this is, on the non-zero coefficients from the regression problems \eqref{eq:modified-Cholesky-decomposition}. Consequently, if we want to force a particular structure on $\BEST^{-1}$ some of the coefficients in \eqref{eq:modified-Cholesky-decomposition} must be set to zero. Thus, we can condition the predecessors of a particular model component to be inside the scope of some radius $\ra$. This will depend on the manner how the model components are labeled. In practice, row-major and column-major formats are commonly used in the context of data assimilation but, other formats can be used in order to exploit particular features of model discretizations and/or dynamics. For instance, making use of row-major format, consider we want to compute the corresponding set of coefficients for the grid point 6 in figure \ref{fig:predecessors} for $\ra=1$. The local box surrounding the grid point 6 provides the model components inside the scope of $\ra$. Readily, the predecessors of 6 are the model components labeled from 1 to 5 according to the labelling system utilized.
In general, the analysis increments of the EnKF reads:
\begin{eqnarray}
\label{eq:analysis-ensemble}
\X^a = \X^b + \DX^a \in \Re^{\Nstate \times \Nens} ,\,
\end{eqnarray}
where $\DX^a$ is known as the analysis increment. According to the primal formulation of the EnKF, $\BEST^{-1}$ is used in order to compute the analysis correction:
\begin{eqnarray}
\label{eq:primal-increment}
\displaystyle 
\DX &=& \lb \BEST^{-1}  + \H^T \cdot \R^{-1} \cdot \H \rb^{-1} \cdot \H^T \cdot \R^{-1} \cdot  \lb \Y^s-\H \cdot \X^b \rb \in \Re^{\Nstate \times \Nens}
\end{eqnarray}
while, in the dual formulation $\BEST$ is implicitly used:
\begin{eqnarray}
\label{eq:dual-increment}
\displaystyle
\DX &=& \X \cdot \V^T \cdot \lb\R + \V \cdot \V^T \rb^{-1} \cdot \lb \Y^s - \H \cdot \X^b \rb \in \Re^{\Nstate \times \Nens} , \,
\end{eqnarray}
where
\begin{eqnarray}
\label{eq:trianlgular-linear-system}
\displaystyle \T \cdot \X = \D^{1/2} \in \Re^{\Nstate \times \Nstate} \,,
\end{eqnarray}
$\Y^s \in \Re^{\Nobs \times \Nens}$ is the matrix of perturbed observation with data error distribution $\Nor \lp \zero_{\Nobs },\, \R\rp$, and $\V = \H \cdot \X \in \Re^{\Nobs \times \Nstate}$. The primal approach can be employed making use of iterative solvers in order to solve the implicit linear system in \eqref{eq:primal-increment}. On the other hand, the dual approach relies most of its computation in the solution of the unitary triangular linear system in \eqref{eq:trianlgular-linear-system}. In general, there are good linear solvers in the current literature, some of them well-known and used in operational data assimilation such as the case of LAPACK \cite{LAPACKLIB} and CuBLAS \cite{BLASLIB}. Compact representation of matrices can be used as well in order to exploit the structures of $\BEST^{-1}$ and $\T$ in terms of memory allocation.

\section{Proposed parallel implementation of the ensemble Kalman filter based on modified Cholesky decomposition}
\label{sec:parallel-implementation}

\begin{figure}[htp]
\centering
\begin{subfigure}{0.45\textwidth}
\centering
\includegraphics[width=0.9\textwidth]{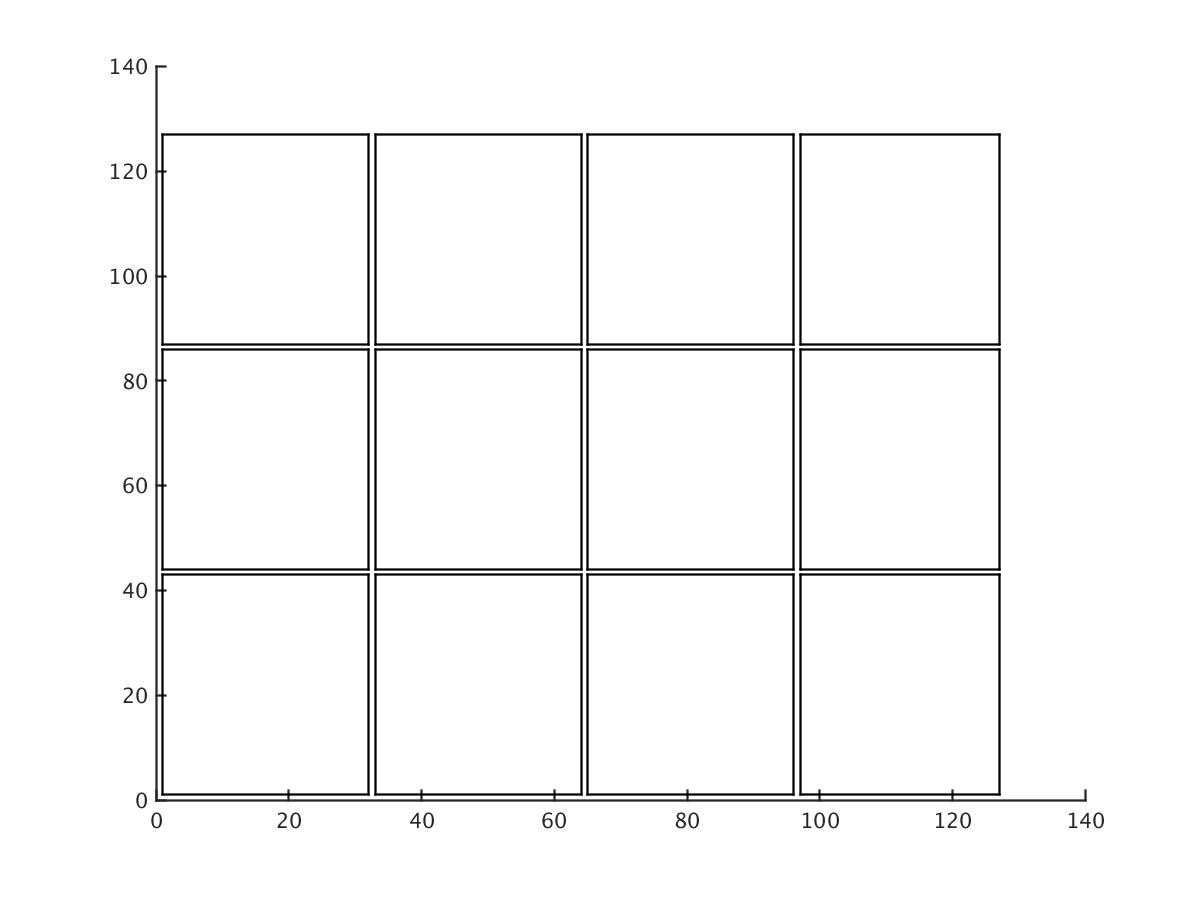}
\caption{Sub-domains 12}
\label{fig:SD12}
\end{subfigure}%
\begin{subfigure}{0.45\textwidth}
\centering
\includegraphics[width=0.9\textwidth]{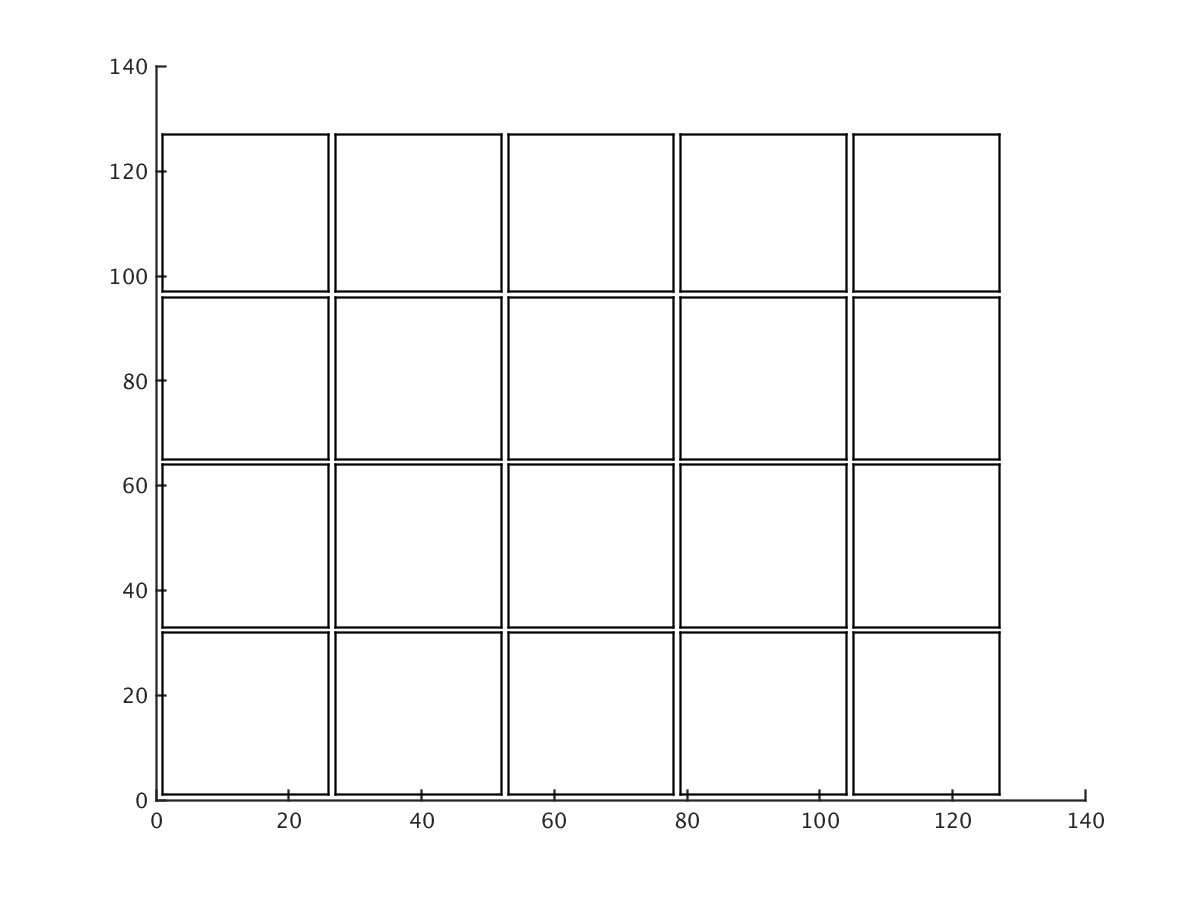}
\caption{Sub-domains 20}
\label{fig:SD20}
\end{subfigure}

\begin{subfigure}{0.45\textwidth}
\centering
\includegraphics[width=0.9\textwidth]{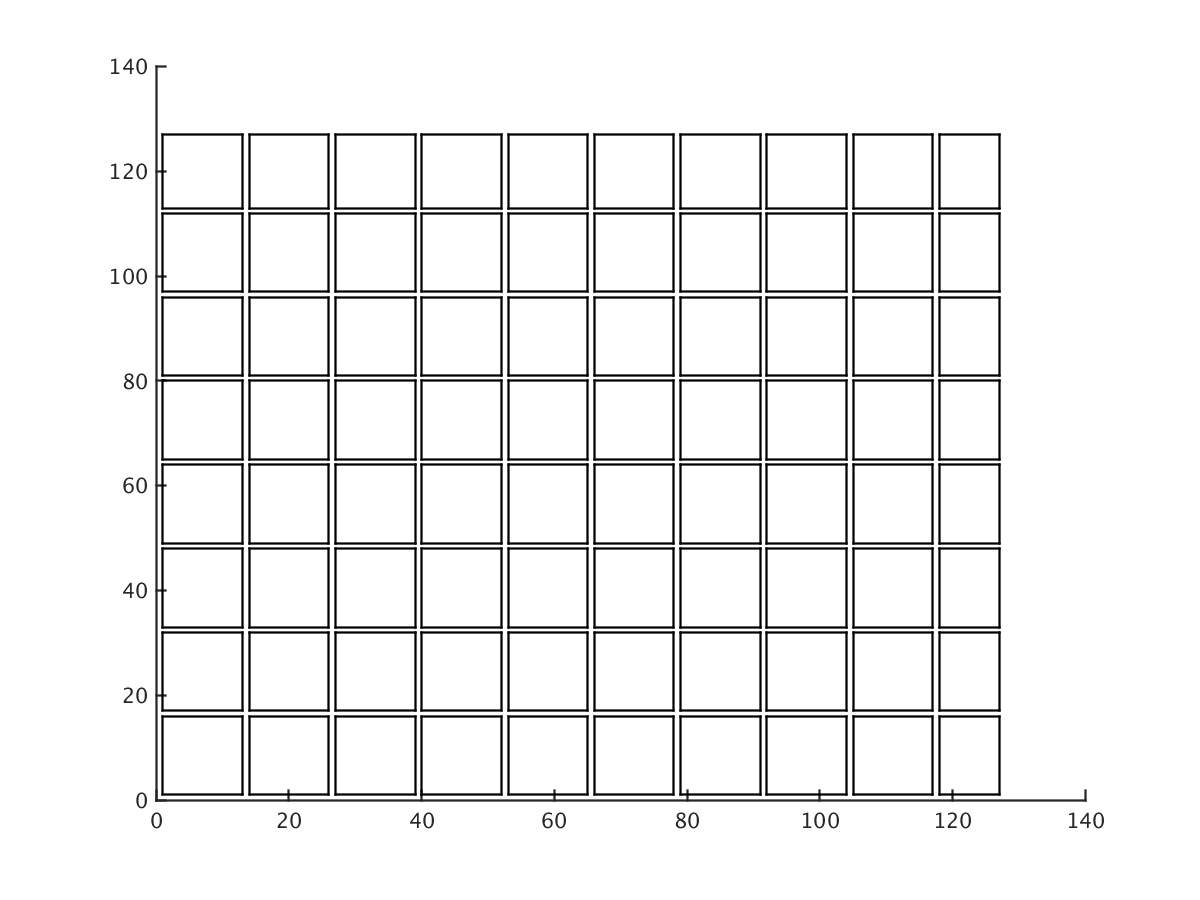}
\caption{Sub-domains 80}
\label{fig:SD80}
\end{subfigure}%
\begin{subfigure}{0.45\textwidth}
\centering
\includegraphics[width=0.9\textwidth]{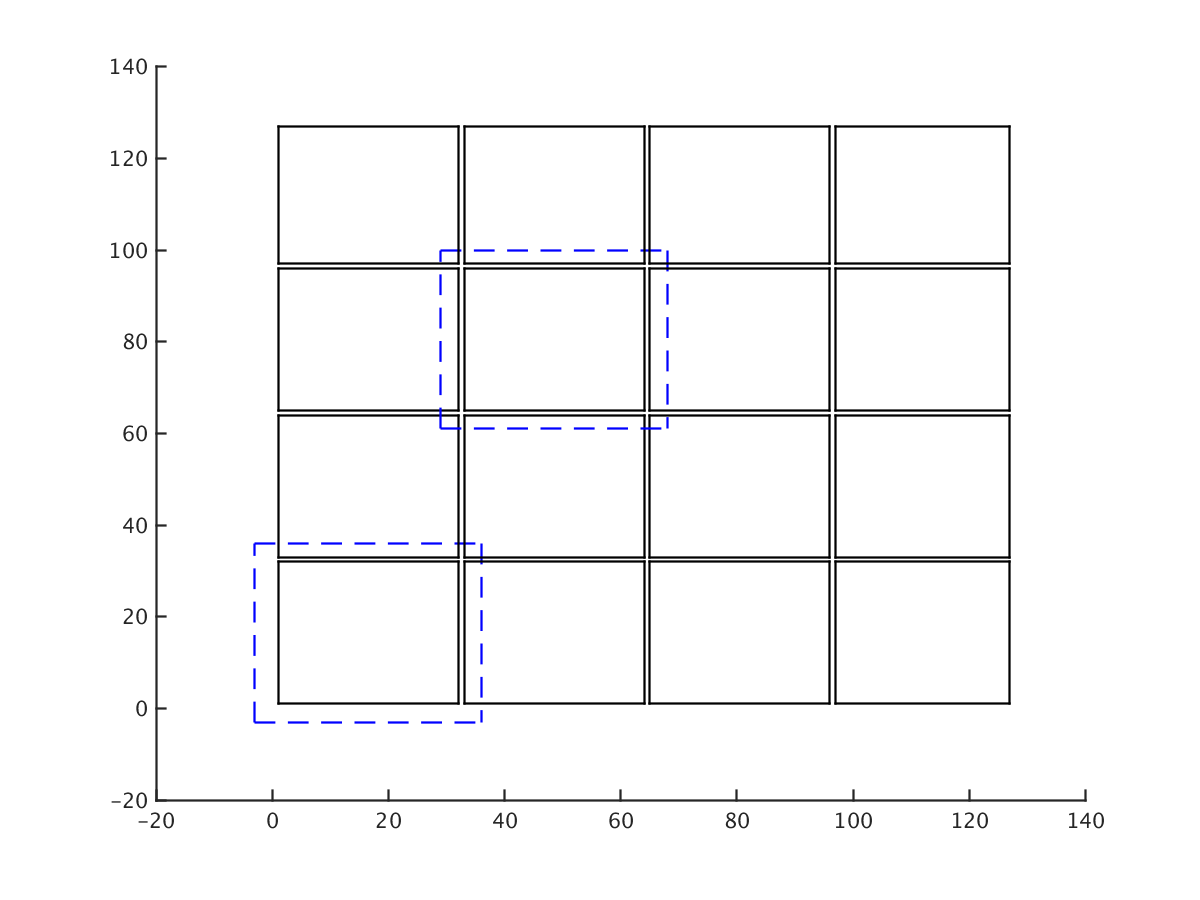}
\caption{Sub-domains 16}
\label{fig:SD16-boundaries}
\end{subfigure}
\caption{Domain decomposition for different number of sub-domains. }
\end{figure}

We consider the use of domain decomposition in order to reduce the dimension of the data assimilation problem. To start, the domain is split according to a given number of sub-domains. Typically, the number of sub-domains matches the number of threads/processors involved in the assimilation process. In figures \ref{fig:SD12}, \ref{fig:SD20} and \ref{fig:SD80} the domain is decomposed in 12, 20 and 80 equitable sub-domains, respectively. With no loose of generality, consider the number of sub-domains $\Delta$ to be a multiple of $\Nstate$. The total number of model components at each sub-domain is $\Nstate / \Delta$ but, in order to estimate $\BEST^{-1}$, boundary information is needed which adds $(2 \cdot \ra+1)^2$ model grid points to the procedure of background covariance matrix estimation. For instance, figure \ref{fig:SD16-boundaries} shows a domain decomposed in 16 sub-domains, the blue dashed squares denote boundary information for two particular sub-domains. 
If we consider $\Delta$ sub-domains, at the $k$-th sub-domain, for $1 \le k \le \Delta$, the analysis reads:
\begin{eqnarray}
\label{eq:analysis-subdomain}
\X^a_{[k]} &=& \X^b_{[k]}+ \BEST_{[k]} \cdot \H^T_{[k]} \cdot \lb \R_{[k]}+\H_{[k]} \cdot \BEST_{[k]} \cdot \H_{[k]}^T \rb  \\ \nonumber
& \cdot & \lb \Y^s_{[k]} - \H_{[k]} \cdot \X^b_{[k]}\rb \in \Re^{\Nstate_{sd} \times \Nens} ,\,
\end{eqnarray}
where $\Nstate_{sd}  =  \Nstate / \Delta +  (2 \cdot \ra+1)^2$, and at sub-domain $k$: $\X^b_{[k]}$ are the model components, $\H_{[k]} \in \Re^{\Nobs_{sd} \times \Nstate_{sd}}$ is the linear observational operator, $\Nobs_{sd}$ is the number of observed components in the sub-domain, $\Y^s_{[k]} \in \Re^{\Nobs_{sd} \times \Nens}$ is the sub-set of perturbed observations, $\BEST^{-1}_{[k]} \in \Re^{\Nstate_{sd} \times \Nstate_{sd}}$ is the local inverse estimation of the background error covariance matrix and $\R_{[k]} \in \Re^{\Nobs_{sd} \times \Nobs_{sd}}$ is the local data-error covariance information. Thus, for all $1 \le k \le \Delta$, the analysis sub-domains \eqref{eq:analysis-subdomain} are computed, the $(2 \cdot \ra + 1)^2$ boundary points are discarded and then, $\Nstate / \Delta$ analysis points are mapped back onto the global domain. Readily, the dual approach can be used as well. One desired property of the proposed EnKF implementation is that boundary information is not exchanged during the assimilation process, each sub-domain works independently in the estimation of $\BEST^{-1}_{[k]}$ and posterior assimilation of $\Y^s_{[k]}$. In the Algorithm \ref{alg:PAR-EnKF-MC}, the parallel ensemble Kalman filter based on modified Cholesky decomposition is detailed. The analysis step of this method is shown in the Algorithm \ref{alg:perform-assimilation} wherein, the model state is divided according to the number of sub-domains $\Delta$ and then, in parallel, information of the background ensemble, the observed components, the observation operator, the estimated data error correlations at each sub-domain are utilized in order to perform the local assimilations. The analysis sub-domains are then merged into the global analysis state as can be seen in line \ref{alg:merging-into-analysis} of the Algorithm \ref{alg:perform-assimilation}. Atomicity is not needed for this operation since analysis sub-domains do not intersect owing to all information concerning to boundaries is discarded after the assimilation step. The local assimilation process is detailed in the Algorithm \eqref{alg:perform-local-assimilation}.
\begin{algorithm}[H]
\caption{Parallel ensemble Kalman filter based on modified Cholesky decomposition (PAR-EnKF-MC)}\label{alg:PAR-EnKF-MC}
\begin{algorithmic}[1]
\Require Initial background ensemble $\X^b = \lb \x^{b[1]},\,\x^{b[2]},\,\ldots,\,\x^{b[\Nens]}\rb \in \Re^{\Nstate \times \Nens}$.
\Ensure Analysis ensemble at each assimilation time.
\While{There are observations to be assimilated}
\State Retrieve $\y$.
\State $\Y^s \gets {\bf create\_perturbed\_observations}(\y,\R)$
\State $\X^{a} \gets {\bf perform\_assimilation}(\X^b,\, \Y^s,\, \R, \,\H)$ \Comment{Parallel analysis step}
\For{${\bf all}\,k \gets 1 \to \Nens$} \Comment{Parallel forecast step}
\State $\x^{b[k]} \gets \M_{t_{previous}\rightarrow t_{current}} (\x^{a[k]})$
\EndFor
\EndWhile
\end{algorithmic}
\end{algorithm}
\begin{algorithm}[H]
\caption{Assimilation step for the PAR-EnKF-MC}\label{alg:perform-assimilation}
\begin{algorithmic}[1]
\Require Background ensemble $\X^b \in \Re^{\Nstate \times \Nens}$, perturbed observations $\Y^s \in \Re^{\Nobs \times \Nens}$, linearized observation operator $\H \in \Re^{\Nobs \times \Nens}$, estimated data error covariance matrix $\R \in \Re^{\Nobs \times \Nobs}$.
\Ensure Analysis ensemble $\X^a \in \Re^{\Nstate \times \Nens}$.
\Procedure{perfofm\_assimilation}{$\X^b$, $\Y^s$, $\R$, $\H$}\Comment{Ensemble members are stored columnwise}
\State Decompose the model states $\X^{b}$ into $\Delta$ sub-domains
\For{${\bf all}\,k \gets 1 \to \Delta$}
\State $\X^b_{[k]} \gets {\bf components\_from\_domain\_k}(\X^b,\,k)$
\State $\H_{[k]} \gets {\bf components\_from\_domain\_k}(\H,\,k)$
\State $\Y^s_{[k]} \gets {\bf components\_from\_domain\_k}(\Y^s,\,k)$
\State $\R_{[k]} \gets {\bf components\_from\_domain\_k}(\R,\,k)$
\State $\X^a_{[k]} \gets {\bf perform\_local\_assimilation}(\X^b,\, \Y^s,\, \R,\, \H)$
\State $\X^a \gets {\bf build\_analysis\_state}(\X^a,\,\X^a_{[k]},\,k)$ \label{alg:merging-into-analysis}
\EndFor
\State \textbf{return} $\X^a$\Comment{The analysis ensemble is $\X^a$.}
\EndProcedure%
\end{algorithmic}
\end{algorithm}
\begin{algorithm}
\caption{Local assimilation method}\label{alg:perform-local-assimilation}
\begin{algorithmic}[1]
\Require Local background ensemble $\X^b_{l} \in \Re^{\Nstate_{sd} \times \Nens}$, local perturbed observations $\Y^s_l \in \Re^{\Nobs_{sd} \times \Nens}$, local linearized observation operator $\H_l \in \Re^{\Nobs_{sd} \times \Nens}$, local estimated data error covariance matrix $\R_l \in \Re^{\Nobs_{sd} \times \Nobs}$.
\Ensure Analysis ensemble $\X^a_l \in \Re^{\Nstate_{sd} \times \Nens}$.
\Procedure{perform\_local\_assimilation}{$\X^b_l$, $\Y^s_l$, $\R_l$, $\H_l$}\Comment{Ensemble members are stored columnwise}
\State Estimate $\BE^{-1}_l$ based on the samples $\X^b_l$.
\State Perform the assimilation,
\begin{eqnarray*}
\displaystyle
\X^a_l \gets \X^b_l+\lb \BE^{-1}_l+\H_l^T \cdot \R^{-1}_l \cdot \H_l \rb^{-1} \cdot \lb \Y^s_l-\H_l \cdot \X^b_l \rb
\end{eqnarray*}
\State \textbf{return} $\X^a_l$\Comment{The local analysis ensemble is $\X^a$.}
\EndProcedure%
\end{algorithmic}
\end{algorithm}

We are now ready to test our proposed parallel implementation of EnKF based on modified Cholesky decomposition.

\section{Experimental Settings}
\label{sec:experiments}

In this section we study the performance of the proposed parallel ensemble Kalman filter based on modified Cholesky decomposition (PAR-EnKF-MC). The experiments are performed using the atmospheric general circulation model SPEEDY \cite{Speedy1,Speedy2}. SPEEDY is a hydrostatic, spectral coordinate, spectral transform model in the vorticity-divergence form, with semi-implicit treatment of gravity waves. The number of layers in the SPEEDY model is 8 and the T-63 model resolution ($192 \times 96$ grids)  is used for the horizontal space discretization of each layer. Four model variables are part of the assimilation process: the temperature ($K$), the zonal and the meridional wind components ($m/s$), and the specific humidity ($g/kg$). The total number of model components is $\Nstate = 589,824$. The number of ensemble members is $\Nens=94$ for all the scenarios. The model state space is approximately 6,274 times larger than the number of ensemble members ($\Nstate \gg \Nens$). The tests are performed on the super computer Blueridge cluster at the university of  Virginia Tech. BlueRidge is a 408-node Cray CS-300 cluster. Each node is outfitted with two octa-core Intel Sandy Bridge CPUs and 64 GB of memory, for a total of 6528 cores and 27.3 TB of memory systemwide.

Starting with the state of the system  $\x^\textnormal{ref}_{-3}$ at time $t_{-3}$, the model solution $\x^\textnormal{ref}_{-3}$ is propagated in time over one year:
\begin{eqnarray*}
\x^\textnormal{ref}_{-2} = \M_{t_{-3} \rightarrow t_{-2}} \lp \x^\textnormal{ref}_{-3}\rp.
\end{eqnarray*}
The reference solution $\x^\textnormal{ref}_{-2}$ is used to build a perturbed background solution:
\begin{eqnarray}
\label{eqch5:exp-perturbed-background}
\displaystyle
\widehat{\x}^{b}_{-2} = \x^\textnormal{ref}_{-2} + \errobs^{b}_{-2}, \quad  \errobs^{b}_{-2} \sim \Nor \lp \zero_{\Nstate} ,\, \underset{i}{\textnormal{diag}} \left\{ (0.05\, \{\x^\textnormal{ref}_{-2}\}_i)^2 \right\} \rp.
\end{eqnarray}
The perturbed background solution is propagated over another year to obtain the background solution at time $t_{-1}$:
\begin{eqnarray}
\label{eqch5:exp-background-state-1}
\x^{b}_{-1} = \M_{t_{-2} \rightarrow t_{-1}} \lp \widehat{\x}^{b}_{-2}\rp.
\end{eqnarray}
This model propagation attenuates the random noise introduced in \eqref{eqch5:exp-perturbed-background} and makes the background state \eqref{eqch5:exp-background-state-1} consistent with the physics of the SPEEDY model. Then, the background state \eqref{eqch5:exp-background-state-1} is utilized in order to build an ensemble of perturbed background states:
\begin{eqnarray}
\label{eqch5:exp-perturbed-ensemble}
\displaystyle
\widehat{\x}^{b[i]}_{-1}  = \x^{b}_{-1} + \errobs^{b}_{-1},\quad \errobs^{b}_{-1} \sim \Nor \lp \zero_{\Nstate} ,\, \underset{i}{\textnormal{diag}} \left\{ (0.05\, \{\x^{b}_{-1}\}_i)^2 \right\} \rp,
\quad 1 \le i \le \Nens,
\end{eqnarray}
from which, after three months of model propagation, the initial ensemble is obtained at time $t_0$:
\begin{eqnarray*}
\x^{b[i]}_0 = \M_{t_{-1} \rightarrow t_0} \lp \widehat{\x}^{b[i]}_{-1}\rp \,.
\end{eqnarray*}
Again, the model propagation of the perturbed ensemble ensures that the ensemble members are consistent with the physics of the numerical model. 

\begin{figure}[H]
\centering
\begin{subfigure}{0.5\textwidth}
\centering
\includegraphics[width=0.9\textwidth,height=0.5\textwidth]{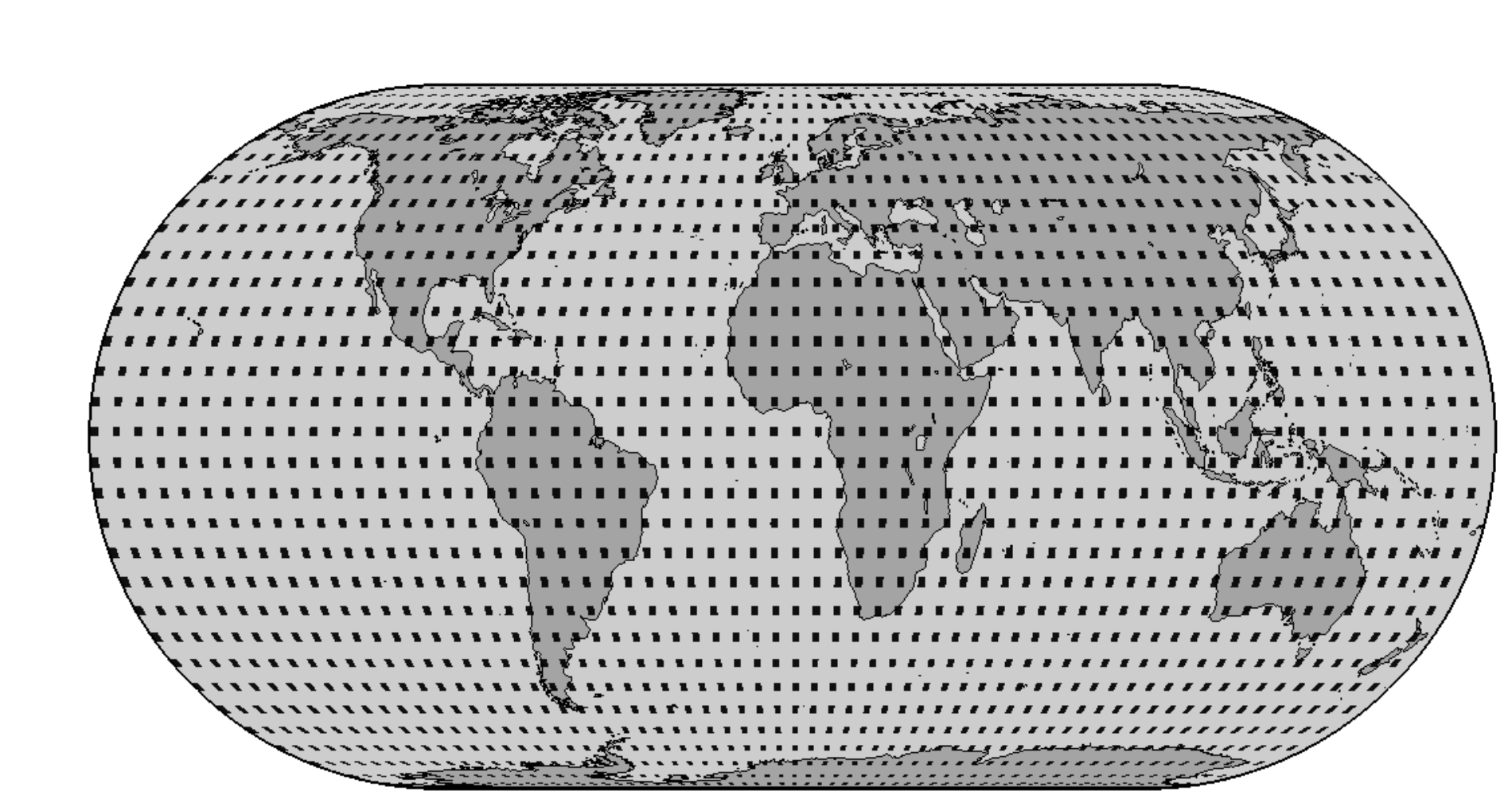}
\caption{$p=12\%$ }
\end{subfigure}%
\begin{subfigure}{0.5\textwidth}
\centering
\includegraphics[width=0.9\textwidth,height=0.5\textwidth]{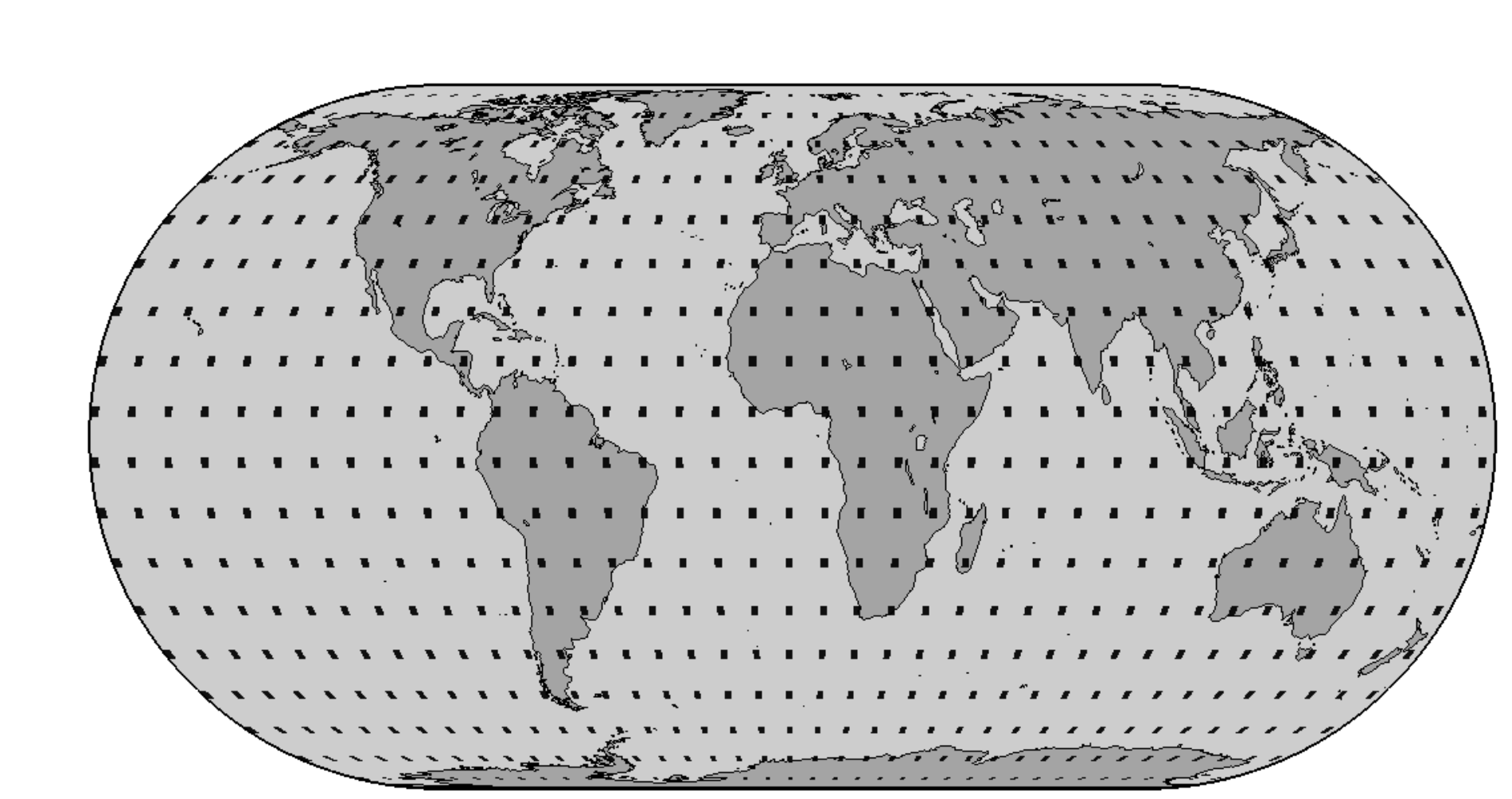}
\caption{$p=4\%$ }
\end{subfigure}%
\caption{Observational networks for different values of $p$. Dark dots denote the location of the observed components. The observed model variables are the zonal and the meridional wind components, the specific humidity, and the temperature.}
\label{fig:exp-observational-grids}
\end{figure}

The experiments are performed over a period of 24 days, where observations are taken every 2 days ($\N=12$). At time $k$ synthetic observations are built as follows:
\begin{eqnarray*}
\y_k = \H_k \cdot \x^\textnormal{ref}_k + \errobs_k, \quad \errobs_k \sim \Nor \lp \zeros_{\Nobs},\, \R_k \rp,\,
\quad \R_k = \textnormal{diag}_{i}\left\{ (0.01\, \{\H_k \, \x^\textnormal{ref}_k\}_i )^2  \right\}.
\end{eqnarray*}
The observation operators $\H_k$ are fixed throughout the time interval. We perform experiments with several operators characterized by different  proportions $p$ of observed components from the model state $\x^\textnormal{ref}_k$ ($\Nobs \approx p \cdot \Nstate$). We consider four different values for $p$: 0.50, 0.12, 0.06 and 0.04 which represent 50\%, 12 \%, 6 \% and 4 \% of the total number of model components, respectively. Some of the observational networks used during the experiments are shown in Figure \ref{fig:exp-observational-grids} with their corresponding percentage of observed components from the model state. 
The analyses of the PAR-EnKF-MC are compared against those obtained making use of the LETKF implementation proposed by Hunt et al in \cite{terasaki2015local,TELA:TELA076,miyoshi2012local} . The analysis accuracy is measured by the root mean square error (RMSE)
\begin{eqnarray}
\label{eqch5:ER-RMSE-formula}
\displaystyle
\text{RMSE}  = \sqrt{\frac{1}{\N} \cdot \sum_{k=1}^\N \lb \x^\textnormal{ref}_k -\x^\textnormal{a}_k \rb^T \cdot \lb \x^\textnormal{ref}_k -\x^\textnormal{a}_k \rb }
\end{eqnarray}
where $\x^\textnormal{ref} \in \Re^{\Nstate \times 1}$ and $\x^\textnormal{a}_{k} \in \Re^{\Nstate \times 1}$ are the reference and the analysis solutions at time $k$, respectively, and $\N$ is the number of assimilation times. 

During the assimilation steps, the data error covariance matrices $\R_k$ are used and therefore, no representativeness errors are involved during the assimilation. The different EnKF implementations are performed making use of FORTRAN and specialized libraries such as BLAS and LAPACK are used in order to perform the algebraic computations. 

\subsection{Influence of the  localization radius on analysis accuracy}
\label{subsec:accuracy-of-proposed-method}
\begin{figure}[H]
\centering
\begin{subfigure}{0.5\textwidth}
\centering
\includegraphics[width=1\textwidth]{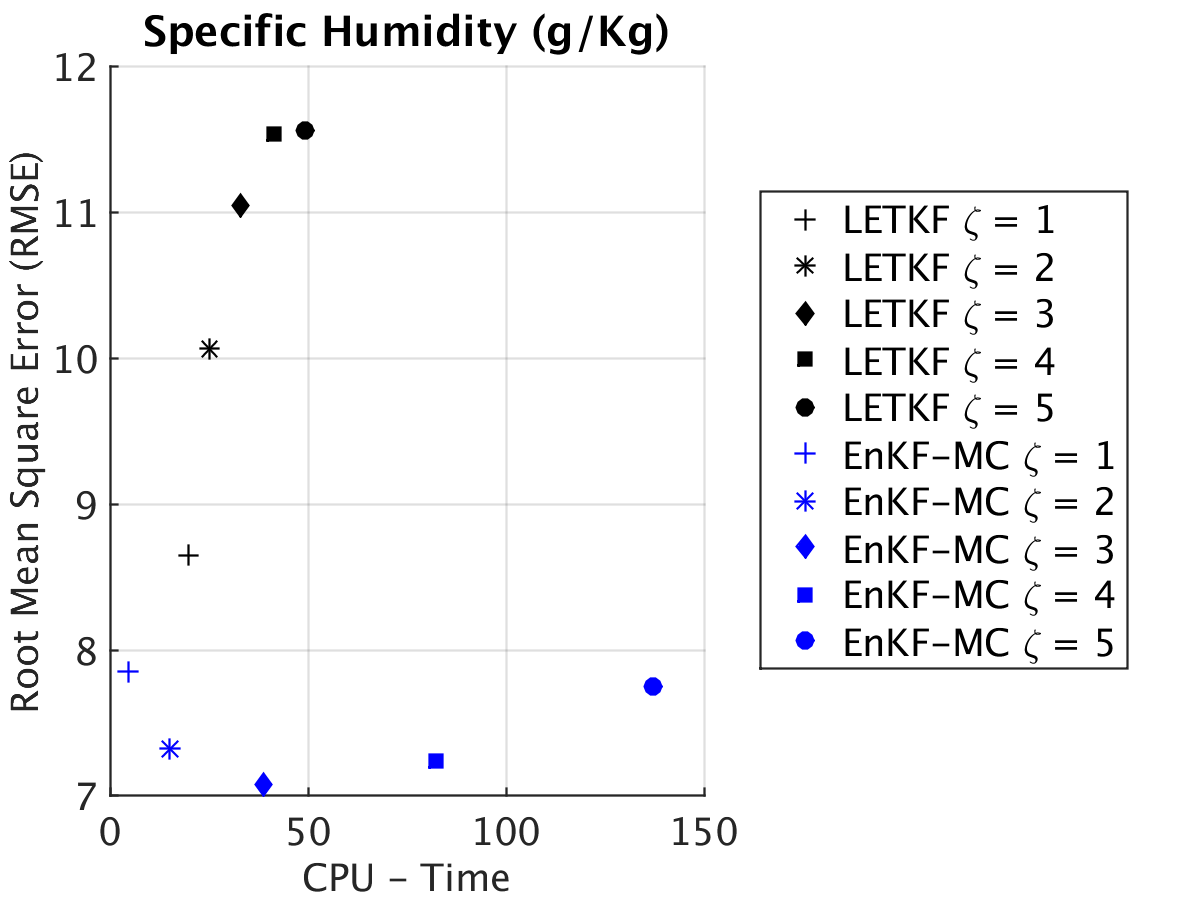}
\caption{$p \sim 50 \%$}
\end{subfigure}%
%
\begin{subfigure}{0.5\textwidth}
\centering
\includegraphics[width=1\textwidth]{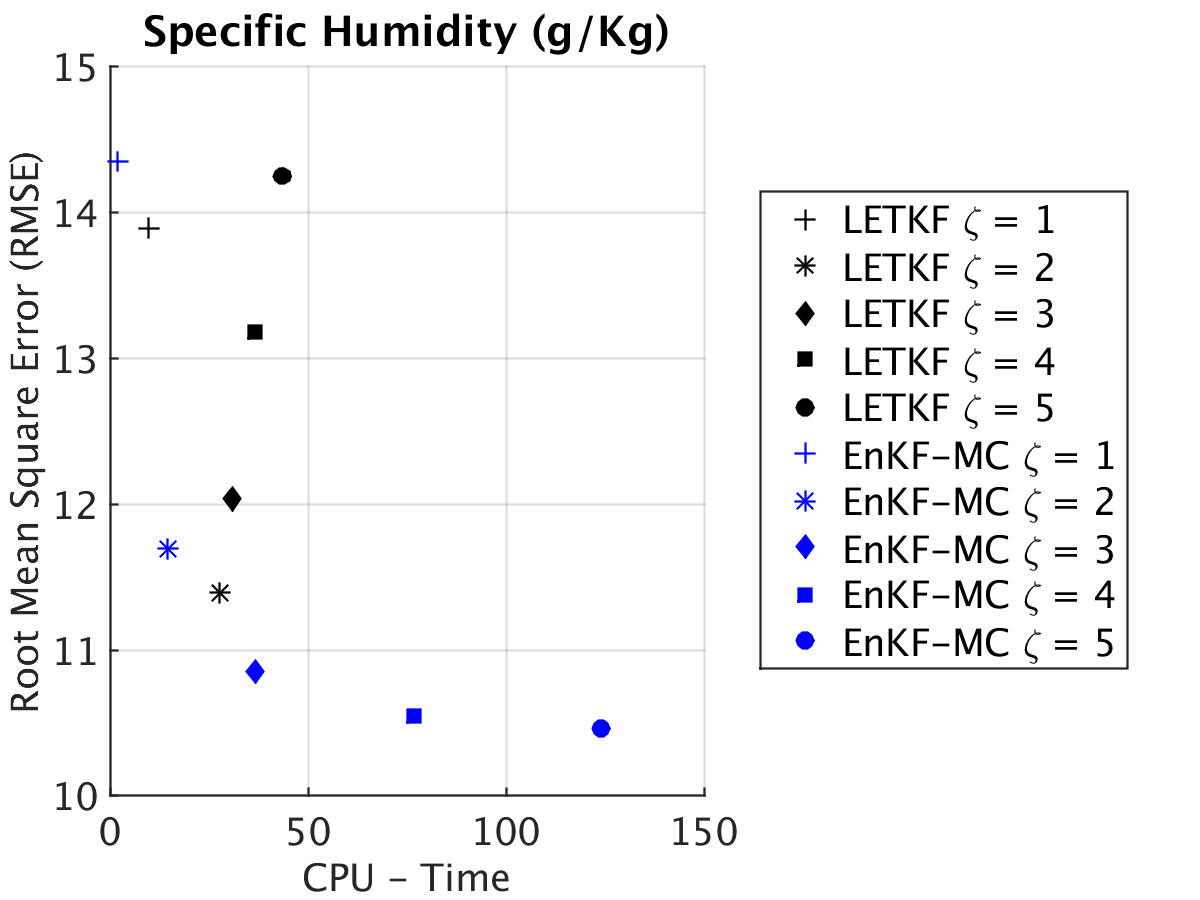}
\caption{$p \sim 4 \%$}
\end{subfigure}%
\caption{Relation between CPU-time (s) and accuracy of the compared EnKF implementations for different radii of influence when the number of computing nodes is 6 (96 processors)}
\label{fig:all-16}
\end{figure}

\begin{figure}[H]
\centering
\begin{subfigure}{0.5\textwidth}
\centering
\includegraphics[width=1\textwidth]{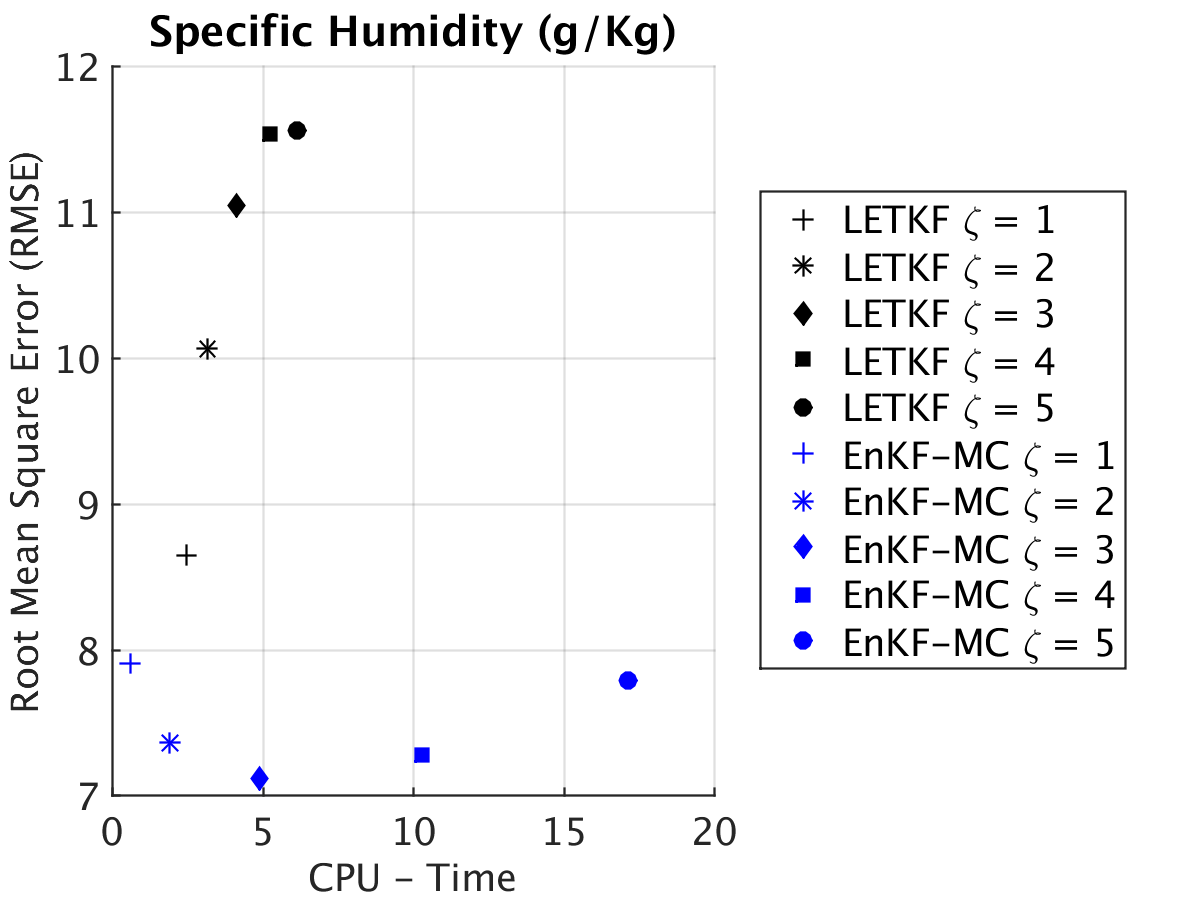}
\caption{$p \sim 50 \%$}
\end{subfigure}%
%
\begin{subfigure}{0.5\textwidth}
\centering
\includegraphics[width=1\textwidth]{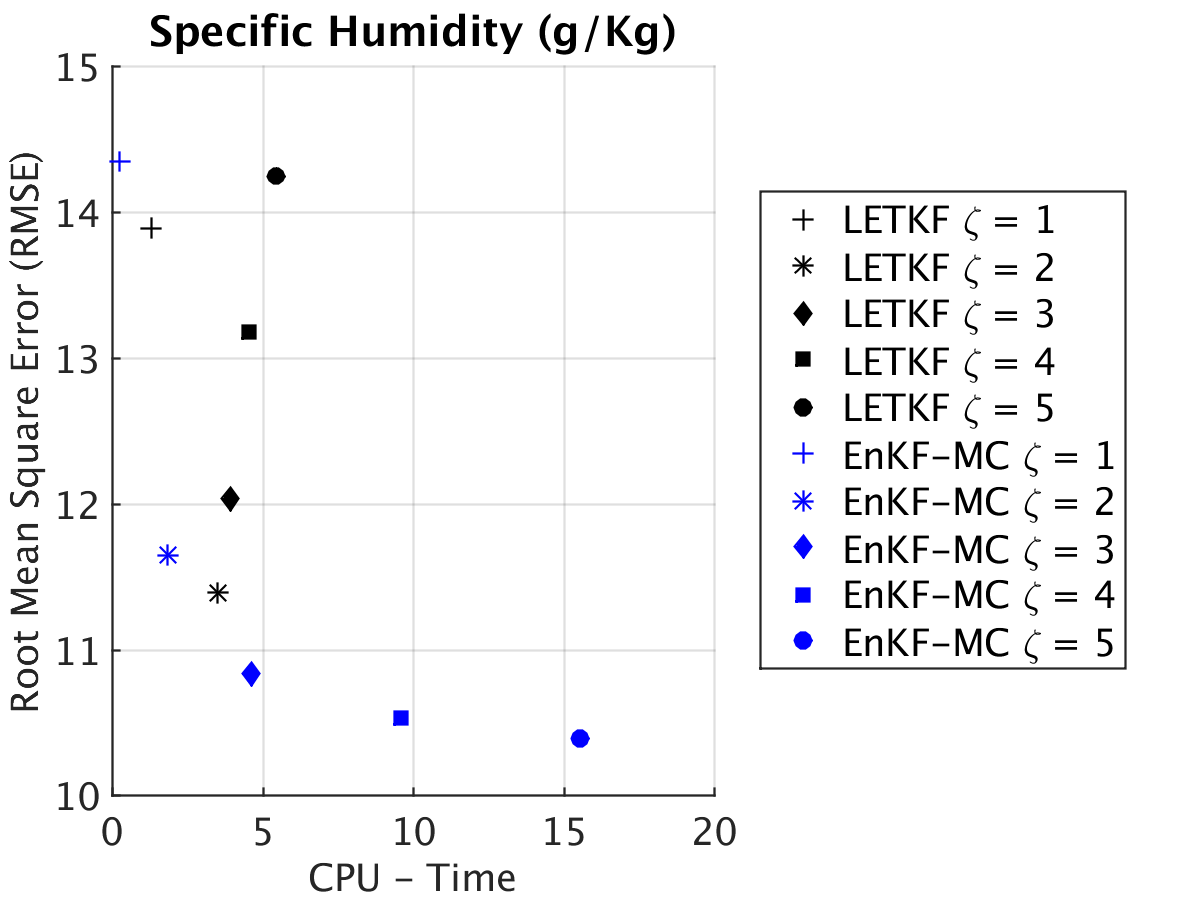}
\caption{$p \sim 4 \%$}
\end{subfigure}%
\caption{Relation between CPU-time (s) and accuracy of the compared EnKF implementations for different radii of influence when the number of computing nodes is 48 (768 processors)}
\label{fig:all-48}
\end{figure}

We study the accuracy of the proposed PAR-EnKF-MC and the LETKF implementations for different radii of influence. The relations between the accuracy of the methods and the radii for 96 and for 768 processors are shown in Figures \ref{fig:all-16} and \ref{fig:all-48}, respectively. The results reveal that the accuracy of the PAR-EnKF-MC formulation can be improved by increasing the radius of influence $\ra$. This implies that the impact of spurious correlations is mitigated when background error correlations are estimated via the modified Cholesky decomposition. However, the larger the radius of influence, the larger the local data assimilation problem to solve. This will demand more computational time which can be mitigated by increasing the number of processors during the assimilation step. On the other hand, in the LETKF context, since background error correlations are estimated based on the empirical moments of the ensemble, spurious correlations affect the analysis when $\ra>2$. Consequently, localization radius sizes beyond this value decreases the performance of the LETKF. 

%
%
\subsection{Computational times for different numbers of processors}
\label{subsec:elapsed-time-per-processor}
\begin{figure}[H]
\centering
\includegraphics[width=0.6\columnwidth]{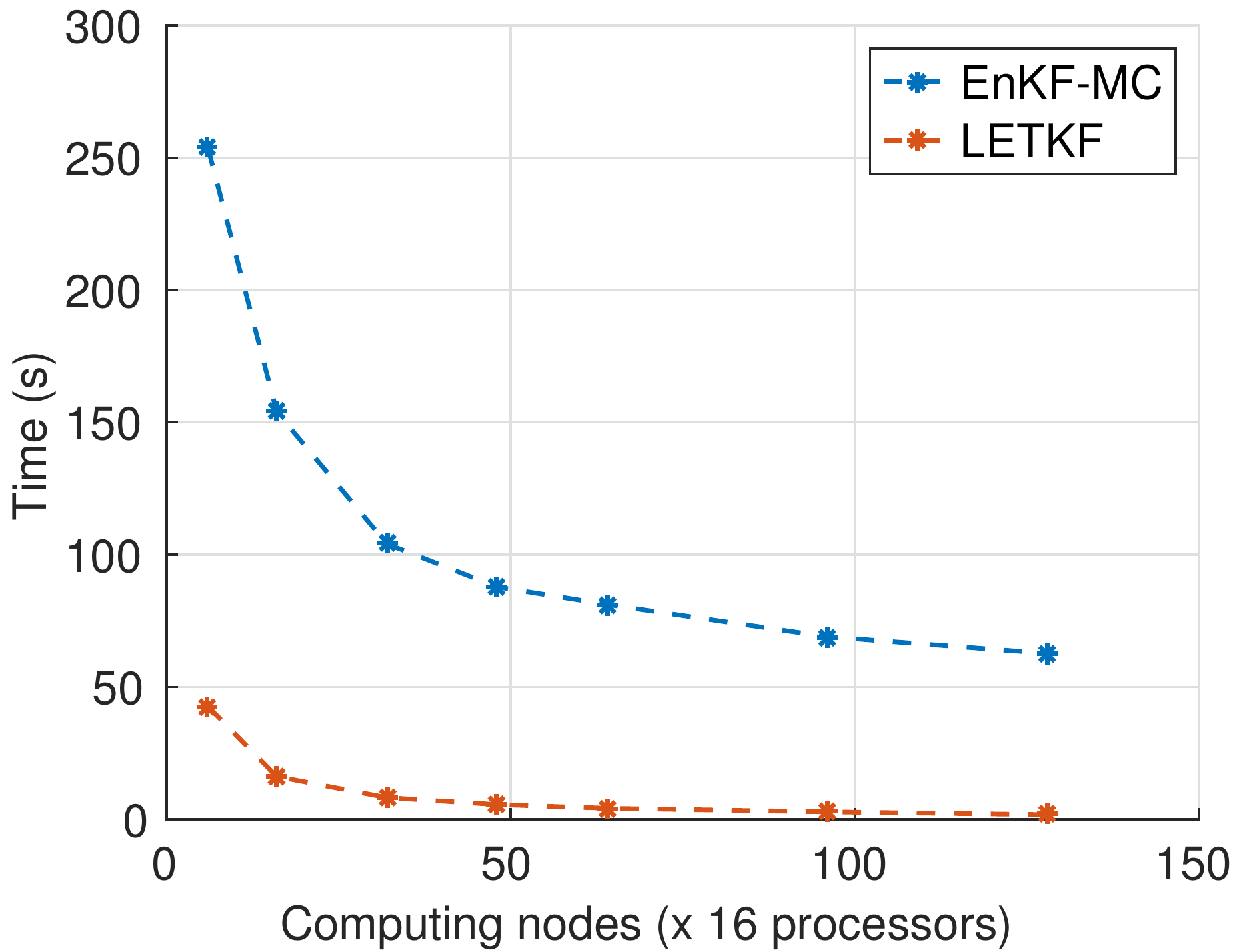}
\caption{Elapsed times of the PAR-EnKF-MC and LETKF for different number of compute nodes ($\times 16$ processors).}
\label{fig:elapsed-times}
\end{figure}

We compare the elapsed times and the accuracy of both implementations when the number of processors (sub-domains) is increased. We vary the number of compute nodes from 6 (96 processors) to 128 (2,048 processors), fix the radius of influence at $\ra=5$, and use an observational network with $p=4\%$. The elapsed times for different numbers of computing nodes for the PAR-EnKF-MC and LETKF are shown in Figure \ref{fig:elapsed-times}. As expected, the elapsed time of the LETKF is smaller than that of PAR-EnKF-MC formulation since no covariance estimation is performed. Nevertheless, the difference between the elapsed times is small (in the order of seconds), while the PAR-EnKF-MC results are more accurate than those obtained by the LETKF.

\subsection{Influence of the number of processors (sub-domains) on accuracy of PAR-EnKF-MC analyses}
\label{subsec:accuracy-PAR-EnKF-MC}
\begin{figure}[htp]
\centering

\includegraphics[width=0.7\textwidth]{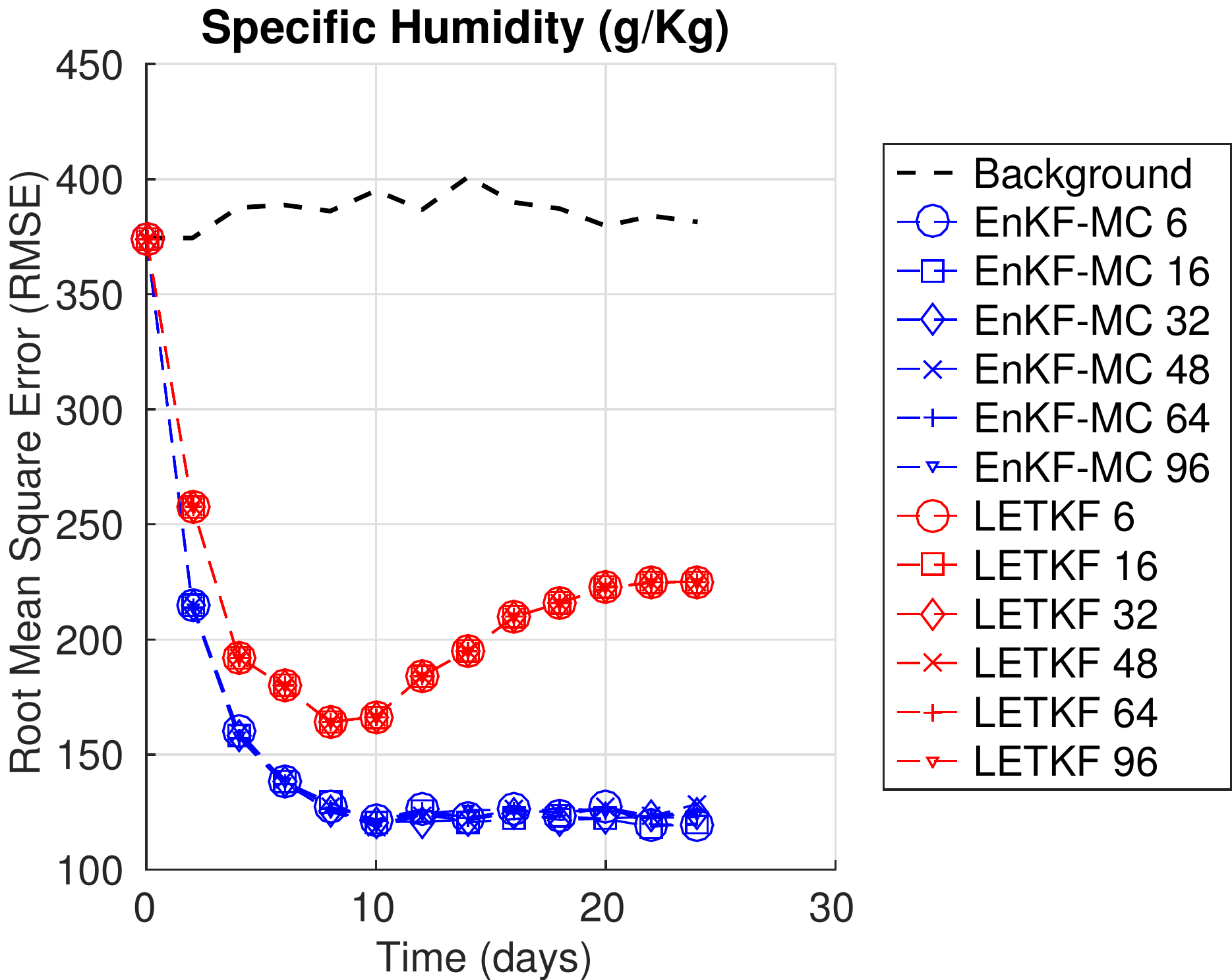}
\caption{RMSE of the LETKF and PAR-EnKF-MC implementations for the specific humidity ($sh$) for different numbers of compute nodes. The number of compute nodes is next to the method name.}
\label{fig:different-variables-processors}
\end{figure}

\begin{figure}[htp]
\centering
\begin{subfigure}{0.45\textwidth}
\centering
\includegraphics[width=0.9\textwidth]{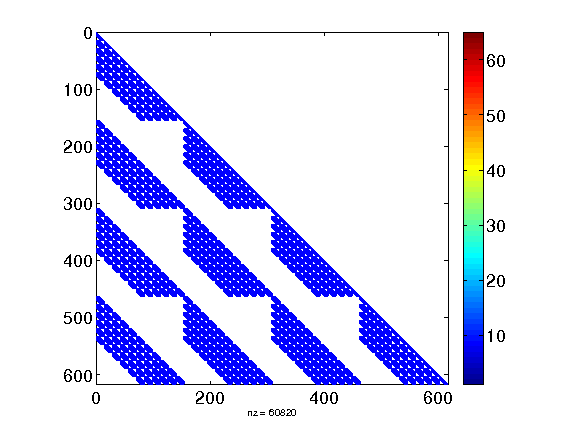}
\caption{Structure of $\T$}
\label{fig:structure-of-T}
\end{subfigure}%
\begin{subfigure}{0.45\textwidth}
\centering
\includegraphics[width=0.9\textwidth]{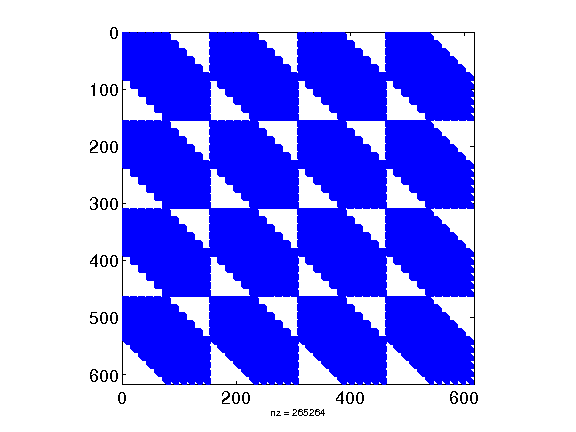}
\caption{Suctructure of $\BEST^{-1}$}
\label{fig:structure-of-B-inv}
\end{subfigure}

\begin{subfigure}{0.45\textwidth}
\centering
\includegraphics[width=0.9\textwidth]{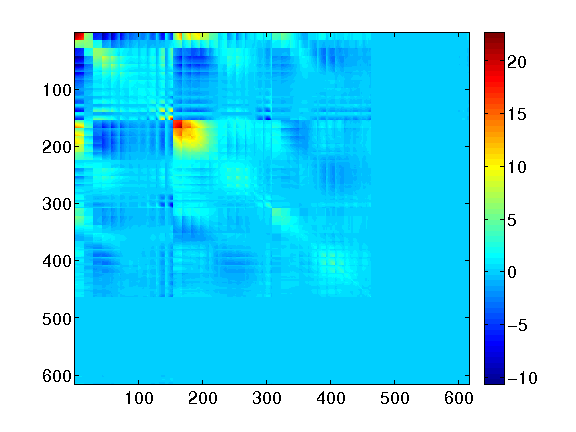}
\caption{$\BEST$}
\label{fig:BEST}
\end{subfigure}%
\begin{subfigure}{0.45\textwidth}
\centering
\includegraphics[width=0.9\textwidth]{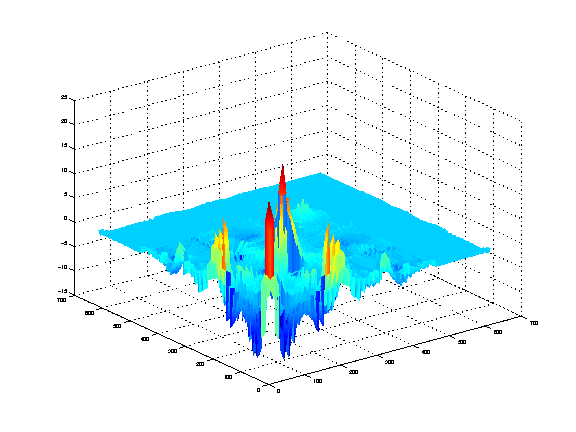}
\caption{Surface of $\BEST$}
\label{fig:BEST-waves}
\end{subfigure}
\caption{Structures of $\T$ and $\BEST^{-1}$ for a radius of influence of $r=5$. The contourf and surface of $\BEST$ are shown as well. The vector state reads $\x = \lb u,\, v,\, T,\, sh \rb^T$.}
\label{fig:local-structure}
\end{figure}

An important concern to address in the PAR-EnKF-MC formulation is how its accuracy is impacted when the number of processors (sub-domains) is increased. As we mentioned before, the model domain is decomposed in order to speedup computations but not for increasing the accuracy of the method (i.e., the impact of spurious correlations can be small for small sub-domain sizes) Two main reasons are that we have a well-conditioned estimated of $\B^{-1}$ and even more, the conditional independence of model components makes the sub-domain size to have no impact in the accuracy of the PAR-EnKF-MC. As can be seen in figure \ref{fig:different-variables-processors}, for the specific humidity variable and values of $\ra$ and $p$, the PAR-EnKF-MC provides almost the same accurate results among all configurations. The small variations in the RMSE values of the PAR-EnKF-MC obey to the synthetic data built at different processors during the assimilation step. For instance, the random number generators used in the experiments depends on the processors id and therefore, the exact synthetic data is not replicated when the number of processors is changed. In the LETKF context we obtain the exact same results for all configurations since it is a deterministic filter and even more, the assimilation is performed for each grid point in the sub-domain.
Lastly, figure \ref{fig:local-structure} shows an estimate of a local inverse background error covariance matrix for some sub-domain. Figure \ref{fig:structure-of-T} shows the non-zero coefficients in that particular sub-domain, figure \ref{fig:structure-of-B-inv} reflects the structure of $\BEST^{-1}$ based on $\T$. Figures \ref{fig:BEST} and \ref{fig:BEST-waves} show the estimated background error covariance matrix $\BEST$ from two different perspectives. As is expected, the correlations are dissipated in space but, they still quite large as can be seen in figure \ref{fig:BEST-waves}. Intuitively, when the sub-domain size is small, high correlations are present between model components owing to their proximity. On the other hand, when the sub-domain size is large, more disipation is expected on the correlation waves of $\BEST$.

\section{Future Work}
\label{sec:future-work}
We think there is an opportunity to exploit even more high performance computing tools in the context of PAR-EnKF-MC. Here, most of the computational time is spent in the estimation of the coefficients in \eqref{eq:modified-Cholesky-decomposition}. The approximation of those coefficients is performed making use of the singular value decomposition (SVD) SVD implementations are highly proposed in the context of accelerating devices such as Many Core Intel (MIC) \cite{Huang13042015} and the Compute Unified Device Architecture (CUDA) \cite{Lahbar5161058}. Since the analysis corrections are computed at each sub-domain independently, each processor (sub-domain) can submit to a given device the information needed in order to solve the linear regression problem \eqref{eq:modified-Cholesky-decomposition}. Once the solution is computed, the device returns the coefficients to the processor which assembles the received information in $\T$. Generally speaking the process is as follows:
\begin{itemize}
\item The domain is split according to $\Delta$ processors (sub-domains)
\item At each sub-domain a local inverse estimation of the background error covariance matrix is computed:
\begin{itemize}
\item Submit the vectors $\x^{[.]}$ to the assigned device in order to compute the weights in the linear regression \eqref{eq:modified-Cholesky-decomposition}.
\item In the device, compute the coefficients making use of SVD.
\item The subdomain receives the coefficients from the device.
\end{itemize}
\item The non-zero coefficients are placed in their respective positions in $\T$.
\item Continue until the coefficients for all local components have been computed.
\item Perform the local assimilation.
\end{itemize}

\section{Conclusions}
\label{sec:conclusions}

An efficient and parallel implementation of the ensemble Kalman filter based on a modified Cholesky decomposition is proposed. The method exploits the conditional independence of model components in order to obtain sparse estimators of $\B^{-1}$ via the modified Cholesky decomposition. High performance computing can be used in order to speedup the assimilation process: the global domain is decomposed according to the number of processors (sub-domains), at each sub-domain a local estimator of the inverse background error covariance matrix is computed and the local assimilation process is carried out. Each sub-domain is then mapped back onto the global domain where then, the global analysis is obtained. The proposed EnKF implementation is compared against the well-known local ensemble transform Kalman filter (LETKF) making use of the Atmospheric General Circulation Model (SPEEDY) with the T-63 resolution in the super computer cluster Blueridge at Virginia Tech. The number of processors is ranged from 96 to 2,048. The accuracy of the proposed EnKF outperforms that of the LETKF. Even more, the computational time of the proposed implementation differs in seconds of the parallel LETKF method in which no covariance estimation is performed. Finally, for the largest number of processors, the proposed method is 400 times faster than its serial theoretical implementation.


\section{Acknowledgments}

This work was supported in part by awards
NSF CCF --1218454,
AFOSR FA9550--12--1--0293--DEF,
and by the Computational Science Laboratory at Virginia Tech.

\end{document}